\nopagenumbers
\magnification=1200
\baselineskip=14pt plus 0.1pt
\pageno=1

\font\tengothic=eufm10
\font\sevengothic=eufm7
\font\fivegothic=eufm5
\newfam\Gothic
\textfont\Gothic=\tengothic
\scriptfont\Gothic=\sevengothic
\scriptscriptfont\Gothic=\fivegothic
\def\go{\fam\Gothic\tengothic}

\font\tenrus=wncyr10
\font\sevenrus=wncyr10 scaled 700
\newfam\Rus
\textfont\Rus=\tenrus
\scriptfont\Rus=\sevenrus

\def\i1{\accent'044i}
\def\e1{\accent'040e}
\def\l1{l{}p1}

\def\mat #1,#2,#3,#4,{\left({#1\atop #3}{#2\atop #4}\right)}
\def\frac #1,#2,{{{#1}\over {#2}}}
\def\bra#1,{{\left\lbrace {#1}\right\rbrace}}

\def\nd{\not\hskip2.2pt\mid}

\def\L{{\cal L}}

\def\det{{\rm det}}


\def\GL{{\rm GL}}
\def\R{{\bf R}}
\def\C{{\bf C}}
\def\Q{{\bf Q}}

\def\Z{{\bf Z}}

\def\al{\alpha}
\def\w{\omega}
\def\W{\Omega}

\def\la{{\lambda}}
\def\La{{\Lambda}}

\def\O{{\cal O}}
 
\def\De{{\Delta}}

\def\ep{{\varepsilon}}

\def\sqrt{\radical"270370 }

\def\oe{\overline }

\def \c{{\go c}}
\def \e{{\go e}}

\def \i{{\bf i}}

\def \GL{{\rm GL}}

\def  \X{{\cal X}}

\def \Ac{{\cal A}}

\def \ord{{\rm ord}}

\def\ph{\varphi}

\def \Res {{\rm Res}}

\def \Gal{{\rm Gal}}

\def\l1{\langle}

\def\Res{{\rm Res}}

\def\End{{\rm End}}

\def\Spec{{\rm Spec}}

\def\nd{\not\hskip2.2pt\mid}

\def\L{{\cal L}}

\def\det{{\rm det}\,}


\def\GL{{\rm GL}}
\def\R{{\bf R}}

\def\C{{\bf C}}
\def\Q{{\bf Q}}

\def\Z{{\bf Z}}

\def\al{\alpha}
\def\w{\omega}
\def\W{\Omega}

\def\la{{\lambda}}
\def\La{{\Lambda}}

\def\O{{\cal O}}

\def\De{{\Delta}}

\def\ep{{\varepsilon}}

\def\sqrt{\radical"270370 }

\def\oe{\overline }

\def \e{{\go e}}

\def \GL{{\rm GL}}

\def  \X{{\cal X}}

\def \Ac{{\cal A}}

\def \ord{{\rm ord}}

\def\leaderfill{\leaders\hbox to 1em{\hss.\hss}\hfill}

\magnification=\magstep1
\vsize=23.5true cm
\hsize=16truecm
\topskip=1truecm
\headline={\tenrm\hfil\folio\hfil}
\raggedbottom
\abovedisplayskip=3mm
\belowdisplayskip=3mm
\abovedisplayshortskip=0mm
\belowdisplayshortskip=2mm
\normalbaselineskip=12pt
\normalbaselines

\null
\vskip 2.5cm

{\bf
\centerline {GENERALIZED KUMMER CONGRUENCES}
\centerline{AND $p$-ADIC FAMILIES OF MOTIVES}

\medskip
\centerline{\bf ALEXEI A. PANCHISHKIN}
}

\bigskip
\bigskip
{\narrower
{\bf Abstract.} We describe some new general
 constructions of $p$-adic $L$-functions attached to certain 
 arithmetically defined
 complex $L$-functions coming from motives over $\Q$ with coefficiens 
 in a number field $T$, $[T:\Q]<\infty$. 
These constructions 
 are equivalent to proving some
 generalized Kummer
 congruences for critical special values of these complex $L$-functions.\par}

\footnote{}{Research at MSRI is suported in part by NSF grant
 DMS-9022140.}

\bigskip

The paper is based on some talks given by the auther during his visit to MSRI in
February and March 1994.  The purpose of this paper is to describe some new general
 constructions of $p$-adic $L$-functions attached to certain 
 arithmetically defined
 complex $L$-functions coming from motives over $\Q$ with coefficiens 
 in a number field $T$, $[T:\Q]<\infty$. 
These constructions 
 are equivalent to proving some
 generalized Kummer
 congruences for critical special values of these complex $L$-functions.

\medskip
The starting point in the theory of $L$-functions is the expansion of the 
 Riemann zeta-function $\zeta (s)$ into the Euler product:
$$
\zeta (s)= \prod_p (1-p^{-s})^{-1}\ \ =\sum_{n=1}^{\infty
}n^{-s}\ \ \ \  (\rm{Re}(s)>1).
$$  
The set of arguments $s$
 for which $\zeta (s)$ is defined can be extended to all
 $s\in \C, s\not= 1$, and we may regard $\C$ as the group of all
 continuous quasicharacters 
$$
\C=\ {\rm Hom}(\R _+^{\times },\C^{\times }),\ \  \ y\mapsto y^s
$$
 of $\R_+^{\times}$. 
The special values
$\zeta (1-k)$ at negative integers  are rational
numbers:
$$
\zeta (1-k)=-{{B_k}\over{k}}\  ,
$$
 where $B_k$ are
 Bernoulli  numbers, which are defined by the formal power series equality
$$
e^{Bt}=\sum_{n=0}^{\infty }{{{B_n}t^n}\over{n!}}={{t\,e^t} 
\over{e^t-1}}\ ,
$$
 and we know 
 (by the classical Sylvester--Lipschitz theorem)  that
$$ 
 c\in {\bf Z}\ \ {\rm implies}\ \  
 c^k(c^k-1){{B_k}\over{k}}\in \bf Z.
$$

The theory of non-Archimedean zeta-functions originates in the work 
 of Kubota and Leopoldt containing $p$-adic interpolation  
 of these special values. 
Their construction turns out to be equivalent to 
 classical Kummer congruences for the Bernoulli numbers,
 which we recall  here in
 the following form. 
Let $p$ be a fixed prime number, $c>1$ an integer
 prime to $p$. 
Put 
$$
\zeta_{(p)}^{(c)}(-k)=(1-p^k)(1-c^{k+1})\zeta (-k)
$$
 and let $h(x)=\sum_{i=0}^{n}\alpha _i\,x^i\in {\bf Z}_p\lbrack
 x \rbrack $ be a polynomial over the ring ${\bf Z}_p$ of $p$-adic integers
 such that
s$$
x\in {\bf Z}_p \Longrightarrow h(x) \in
p^m {\bf Z}_p.
$$
Then we have that
$$
\sum_{i=0}^{n}\alpha _i\, \zeta_{(p)}^{(c)}(-i)\in p^m {\bf Z}_p.
$$ 
This property expresses the fact that 
 the numbers $\zeta _{(p)}^{(c)}(-k)$ depend continuously on $k$ in the
 $p$-adic sense; it can be  deduced from the known formula for the sum of 
 $k$-th powers:  
$$
S_k(N)=\sum_{n=1}^{N-1}n^k\ \ =\ \ {{1}\over{k+1}}\lbrack
B_{k+1}(N)-B_{k+1}\rbrack
$$ in which  $B_{k}(x)=(x+B)^k=\sum_{i=0}^k{k \choose
 i}B_ix^{k-i}$ denotes the Bernoulli polynomial. 
Indeed, all summands in $S_k(N)$
 depend $p$-adic analytically on $k$, if we restrict ourselves to numbers $n$, 
 prime to $p$, so that the desired congruence follows if we express the numbers
 $\zeta _{(p)}^{(c)}(-k)$ in terms of Bernoulli numbers.

The domain of definition of $p$-adic zeta functions is the $p$-adic
 analytic Lie group
$$X_p={\rm Hom}_{\rm contin}({\bf Z}_p^{\times },{\bf
 C}_p^{\times })
$$ of all continuous $p$-adic characters of the profinite group
$\Z_p^{\times}$,
where ${\bf C}_p={\widehat {\overline {\bf Q}}}_p $ 
 denotes the Tate field (completion of an algebraic
 closure of the $p$-adic field ${\bf Q}_p)$, so that all integers $k$ can be
 regarded  as the characters $x_p^k:\ y\mapsto y^k$. 
The construction of
 Kubota and Leopoldt is equivalent to existence a $p$-adic analytic function
 $\zeta _p:\  X_p \to {\bf C}_p$  with a single pole at the point
 $x=x_p^{-1}$,
 which becomes a bounded holomorphic function on $X_p$
 after multiplication
 by the elementary factor
 $(x_p\,x-1)\ \ (x\in  X_p)$, and is {\sl uniquely
 determined} by the condition 
$$
\zeta _p({x_p}^k)=(1-p^k)\zeta (-k)\ \ \ \ \ \ \
(k\ge 1).
$$
 
This result has a very natural interpretation in framework of the theory of
 non-Archimedean integration (due to Mazur): there exists a $p$-adic measure
$\mu ^{(c)}$ on ${\bf Z}_p^{\times } $ with values in ${\bf Z}_p$ such that
$\int _{{\bf Z}_p^{\times }}x_p^k\,\mu
 ^{(c)}\  = \ \zeta _{(p)}^{(c)}(-k)$. 
Indeed, if we  integrate $h(x)$ over  
 ${\bf Z}_p^{\times } $ we exactly get the above congruence. On the other
 hand, in order to define  a measure $\mu ^{(c)}$ satisfying the above
 condition
 it  suffices  for any continuous function
 $\phi{}:{\bf Z}_p^{\times}\to {\bf Z}_p$
 to define its integral
 $\int _{{\bf Z}_p^{\times }}\phi (x)\,\mu
 ^{(c)}.$  
For this purpose we approximate $\phi (x)$ by a polynomial (for
 which
 the integral is already defined), and then pass to the
 limit.
 
The important feature of the construction is that it equally works 
 for primitive Dirichlet characters $\chi$ modulo a power of $p$: if we fix
 an embedding $i_p:{\overline {\bf Q}}\hookrightarrow
 {\bf C}_p$ then  the character  $\chi{}:({\bf Z}/{\bf Z}_{p^N})^{\times} 
 \to (\overline {\bf Q})^{\times}
  $ becomes an element of the torsion
 subgroup $ X_p^{\rm{tors}}\subset X_p$ and the above equality  also
 holds for the special values $L(-k,\chi )$ of the Dirichlet
 $L$-series
$$
L(s,\chi )= \sum_{n=1}^{\infty}\chi (n)\,n^{-s}=\prod_p (1-\chi
(p)p^{-s})^{-1},
$$
 so that we have 
$$
\zeta _p(\chi x_p^k)=i_p\lbrack 
 (1-\chi (p)p^k) L(-k,
 {\chi })\rbrack\ \ \ \  ( k\ge 1,\ k\in {\bf Z},\ \chi \in 
 X_p^{\rm{tors}}).
$$
The original construction of Kubota and Leopoldt was
 successesfully used by Iwasawa
 for the description of the class groups of cyclotomic fields. 
Since then the class
 of functions admitting
 $p$-adic analogues has gradually extended.

$L$-functions (of complex variable) can be attached as certain
 Euler products to various objects such as  diophantine equations,
 representations of Galois groups, modular forms etc., and they
 play a crucial role in modern number theory. 
Deep interrelations
 between these objects discovered in last decades are based on
 identities for the corresponding $L$-functions which presumably all
 fit into a general concept of the Langlands of $L$-functions associated
 with automorphic representations of a  reductive group $G$ over a
 number field $K$. 
>From this point of  view the study of arithmetic
 properties of these zeta function is becoming especially important.

The major sources of such $L$-functions are: 

\medskip
\noindent
1) {\it  Galois representations } of $G_K=\Gal(\oe K/K)$ for algebraic
 number fields $K$, $r:G_K\to \GL(V)$, $V$ a finite dimensional vector space, 
 and one can attach to $r$ an Euler product due to Artin.

\medskip
\noindent
2) {\it Algebraic varieties $X$} defned over an algebraic number field $K$.
In this case one can attach to $X/K$ its Hasse--Weil zeta function.

\medskip
\noindent
3) {\it Automorphic forms and automorphic representations.}
In the classical case one associates to a modular form $f(z)=
 \sum_{n=0}^\infty a_n\exp (2\pi i nz)$ its Mellin transform
 $L(s,f)=\sum_{n=1}^\infty a_n n^{-s}$.
In general an automorphic form generates an automorphic representations
 in the space of smooth functions over an adelic reductive group, 
 and one can attach an Euler product to it using a decomposition of such a
 representation into a tensor product indexed by prime numbers $p$ and $\infty$.

\medskip
Conjecturally, all the three type of $L$-functions can be related to each
 other
 using a general theory of motives over $\Q$ with coefficiens 
 in a number field $T$, $[T:\Q]<\infty$ (this field coincides with the field
 $\Q(\{a_n\}_{n\ge 1})$ generated by the coefficients of the corresponding 
 $L$-function $L(M,s)=\sum_{n=1}^\infty a_n n^{-s}.$)
For a fixed prime number $p$ one can also attach in many cases to the above
 complex $L$-function 
 a $p$-adic $L$-function. 
These $p$-adic $L$-functions are certain analytic functions in a $p$-adic
 domain obtained by an interpolation procedure of certain special values
 of the corresponding complex analytic $L$-functions.
Their existence is equivalent to certain generalized Kummer
 congruences for the special values. 

\medskip
We describe a general conjecture on the existence of a $p$-adic family $M_P$ of
 motives, parametrised by some dense subset of
 algebraic characters $P$ of a $p$-adic commutative algebraic group (which we
 call group of Hida).
This group can be regarded as (a maximal torus of) the $p$-adic part $G_{M,p}$
of the motivic Galois group $G_M$ of $M$ (the Tannakian group for the tensor
category generated by $M$).  
The important condition of motives $M_P$ of the above family is that they have
the same fixed
 $p$-invariant $h=h_P$ (generalized slope), which is defined  as the difference
 between the Newton
 polygon and the 
 Hodge polygon of a motive at certain point $d^+$ (the dimension of the
 subspace $M^+$ of the Betti realization $M_B$ of $M$).  The corresponding
$p$-adic
 $L$-functions of this family can be unbounded in the "cyclotomic direction" (of
Amice--V\'elu type [Am-Ve])
 but they
 form a family which is conjectually bounded in the "weight direction", that is
 for $P$ parametrized by algebraic characters of $G_{M,p}$.

\medskip
More precisely, the values of the function $P\mapsto L(M_P, 0)$ satisfy 
 generalised Kummer congruences in the followinfg sense:
 for any finite linear combination $\sum_Pb_P\cdot P$ with $b_P\in \C_p$ which
 has the property 
 $\sum_Pb_P\cdot P\equiv 0 (\bmod p^N)$ we have that for some 
 constant $C\not = 0$ the
 corresponding
 linear combination of the normalized $L$-values
$$
C\sum _P b_P c_p(M_P) \cdot \frac  {L_{(p,\infty)}(M_P, 0)}, c_\infty(M_P), 
\equiv 0 (\bmod p^N). 
$$
Here $c_p(M_P)$ and $c_\infty(M_P)$ denote a $p$-adic and a complex period of
 $M_P$ so that the ratio "$\displaystyle \frac c_p(M_P), c_\infty(M_P),$" 
 is uniquely defined, and $L_{(p,\infty)}(M_P, s)$ denotes the above
$L$-function $L(M_P, s)$ normalized by multiplying by a certain canonicaly
defined Deligne's $p$-factor corresponding to a choice of inverse roots 
$\al^{(1)}(p)$, $\dots$, $\al^{(d^+)}(p)\in \C_p$ of  $p$-local
polynomial of $M$ such that
$$
\ord_p(\al^{(1)}(p))\le \ord_p(\al^{(2)}(p))\le \dots \ord_p(\al^{(d)}(p)),
$$ 
$d$ being the common rank of the family $M_P$, $d^+$ the $T$-dimension of the
Deligne'a subspace $M^+$ of $M_B$ (the fixed subspaces of the canonical
involution $\rho$ of $M$ over $T$).

\medskip
Recent examples of such families related to modular forms were constructed by
R.Coleman [CoPBa] who proved the following 

\medskip\noindent
{\bf Theorem. } {\sl Suppose $\al \in \Q$ and $\ep :(\Z/pZ)^\times \to
\C_p^\times$ is a character. 
Then there exists a number $n_0$ which depends on
$p$, $N$ and $\ep$, and $\al$ with the following property:
If $k\in \Z$, $k>\al+1$ and there is a unique normalized cusp form $F$ on
$X_1(Np)$ of weight $k$, character $\ep \w^{-k}$ and slope $\al$ and if
$k^\prime > \al +1$ is an integer congruent to$k$ modulo $p^{n+n_0}$, for any
positive integer $n$, then there exists a unique normalized cusp form $F^\prime$
on $X_1(Np)$ of weight $k^\prime$, character $\ep \w^{-k^\prime}$ and slope
$\al$ ($\w$ denotes the Teichm\"u ller character). 
Moreover his form satisfies
the congruence 
$$
F^\prime(q)\equiv F^\prime(q) (\bmod p^{n+1}).
$$
}

\medskip\noindent
This result can be regarded as a generalization of the work of Hida [HiGal] who
considered the case $\al =0$ and constructed interesting families of
Galois representations of the type  
$$
\rho_p : G_{\Q} \rightarrow \GL_2(\Z_p[[T]]),\ \ G_{\Q} = \Gal(\oe \Q/\Q),
$$
which are non ramified  outside $p$. 
These representations have the following property: if we consider the
homomorphisms 
$$
\Z_p[[T]]\ {\buildrel s_k\over \longrightarrow} \Z_p, \ \ 1 + T \ 
\mapsto (1+p)^{k-1},
$$
then  we obtain a family of Galois representations
$$
\rho_p^{(k)} : G_{\Q} \rightarrow \GL_2(\Z_p),
$$
which is parametrized by  $k \in \Z$, and for $k = 2, 3, \cdots, $ these
representations are equivalent over $\Q_p$  to the $p$-adic representations of
Deligne, attached to modular forms of weight $k$. 
This means that the
representations of Hida are obtained by the $p$-adic interpolation of Deligne's
representations. A geometric interpretation of Hida's  representations was given
by Mazur and Wiles [Maz--W], cf. [Maz].  
For example, for the modular form
$\De$ of weight 12 Hida has constructed a representation 
$$ \rho_{p, \De} :
G_{\Q} \rightarrow \GL_2(\Z_p[[T]]), 
$$ 
as an example of his general theory,
where the prime number $p$ have the property $\tau (p) \not \equiv 0(\bmod p)$
(e.g. $p < 2041, p \not = 2, 3, 5$ and 7).

The boundednes property is the subject of a research
by G.Stevens and of a forthcoming paper of B.Mazur and F.Q.Gouv\^ea.
 
Note that other examples may include Rankin products, Garrett triple products of
 elliptic and Hilbert modular
 forms and standard
 $L$-functions of Siegel modular forms.

\medskip
To describe this conjecture more precisely, let $M$ is a motive over $\Q$ of 
with coefficients in  $T$ i.e.
$$
M_{B},\    M_{DR},\    M_{\lambda },\    I_{\infty},\    I_{\lambda},
$$  
 where 
$M_B$ is the Betti realization of $M $ which is a vector space over $T$ of
dimension $d$ endowed with a $T$-rational
involution $\rho$ ; $M_B = M^+\oplus M^-$ denotes the corresponding
decomposition into the sum of $(+1)$-eigenspaces and $(-1)$-eigenspaces of $\rho$.

\par \medskip
\noindent
$M_{DR}$ is the de Rham realization of $M$, a free $T$-module of rank
$d$, endowed with a decreasing filtration 
$\{F^{i}_{DR}(M)\subset M_{DR}\mid i\in {\bf Z}\}$ of
$T$-modules;
\par\noindent
$M_{\lambda }$ is the $\lambda $-adic realization of $M$ at a finite place 
$\lambda $ of the
coefficient field $T$ (a $T_{\lambda } $-vector space of degree $d$ over 
$T_{\lambda }$, a
completion of $T$ at $\lambda )$ which is a Galois module over $G_{\Q} =
\Gal({}\oe \Q /\Q)$ so that we a have a compatible system of $\lambda $-adic
representations denoted by
\par
$$
r_{M ,  \lambda } = r_{\lambda } :   G_{\Q} \rightarrow  GL(M_{\lambda }).
$$
 Also,
\par
$$
I_{\infty } :   M_B\otimes _T{\bf C}\rightarrow  M_{DR}\otimes
_{ T }{\bf C} 
$$
 is the complex comparison isomorhism of complex vector spaces
$$
I_{\lambda   } :   M_B\otimes _{T}T_{\lambda }\rightarrow  M_{\lambda }
$$
 is the $\lambda $-adic comparison isomorphism of
 $T_{\lambda }$-vector spaces. 
It is
assumed in the notation that the complex vector space 
$M_B\otimes _{{\bf Q}}{\bf C}$ is
decomposed in the Hodge bigraduation
\par
$$
M_{B}\otimes _T{\bf C} = \oplus _{i ,  j}M^{i ,  j}
$$
 in which $\rho (M^{i ,  j}) \subset  M^{j ,  i}$  and 
$$
h (i,  j ) = h (i,  j,   M) = \dim _{{\bf
C}}M^{i ,  j}
$$ are the Hodge
numbers. 
Moreover, \par 
$$
I_{\infty }(\oplus _{i^\prime \ge i}M^{i^\prime  ,  j}) =
F^{i}_{DR}(M)\otimes {\bf C}. 
$$
 Also, $I_{\lambda }$ takes $\rho$ to the $r_{\lambda }$-image of
the Galois automorphism which corresponds to the complex conjugation of ${\bf
C}$.
 We assume that $M$ is pure of weight $w $ (i.e. $i+j = w$).
\par
The $L $-function $L(M,  s)$ of $M$ is defined as the following Euler
product:
$$
L(M,  s) = \prod_{p}L_{{ p}}(M,  { p}^{-s}),
$$  
 extended over all primes ${ p}$  and where
\par
$$
\eqalign{&
L_{{ p}}(M ,  X)^{-1} = \det (1-X\cdot r_{\lambda }(Fr^{-1}_{{ p}})\mid
M^{I_{{ p}}}_{\lambda }) = \cr &
(1-\alpha ^{(1)}({p})X)\cdot (1-\alpha ^{(2)}({p})X)\cdot \dots\cdot
(1-\alpha ^{(d)}({p})X) = \cr &
1+A_{1}({p})X + \dots  + A_{d}({p})X^{d};}
$$
 here  $Fr_{{ p}} \in  G_{\Q}$ is the Frobenius element at
${ p}$, defined modulo conjugation and modulo the inertia subgroup
$I_{{ p}}\subset G_{{p}}\subset G_{\Q}$ of the decomposition group 
$G_{{ p}} $ (of any extension of ${p}$ to
$\oe \Q $). 
We make the standard hypothesis that the coefficients of
$L_{{p}}(M,  X)^{-1}$ belong to $T$, and that they are independent of $\la$
coprime to $p$. 
Therefore we can and we shall regard this polynomial both
over $\C$ and over $\C_p$. 
We shall need the following twist operation: for an
arbitrary motive $M$ over $\Q$ with coefficients in $T$ an integer $m$ and a
Hecke character $\chi $ of finite order one can define the twist 
$N = M(m)(\chi)$ which is again a motive over $\Q$ with the coefficient field
$T(\chi )$ of the same rank $d$ and weight $w$ so that we have 
$$
L(N ,  s) =\prod_{p}L_{{ p}}(M ,  \chi ({ p}) p^{-s-n}).
$$

\medskip
\noindent
{\bf  The group of Hida and the algebra of Iwasawa--Hida.}
Now let us fix a motive $M$ with coefficients in 
$T = \Q (\langle a(n )_{n}\rangle )$ of rang $d$ and of weight $w$, and let
$\End_TM$ denote the endomorphism algebra of $M$ (i.e. the algebra of $T$-linear
endomorphisms of any $M_B$, which commute with the Galois action under the
comparison isomorphisms).
 Let 
$$
\Gal_{p} = \Gal (\Q_{p, \infty}^{ab}/ \Q)
$$ 
 denotes the Galois group of 
 the maximal abelian extension  $\Q_{p, \infty}^{ab}$ of $\Q$ unramified outside
 $p$ and $\infty$. Define $\O_{ T, p} =
\O_T\otimes\Z_p$.

\medskip
\noindent
{\bf Definition.} {\sl The {\it group of Hida } $GH_M = GH_{M, p}$ 
is the following product
$$
GH_M = (\End_TM)^{\times}(\O_{T, p})\times\Gal_p,
$$
where $(\End_TM)^{\times}$ denotes (a maximal torus of) the algebraic $T$-group
of  invertible elements of $\End_TM$ (it is implicitely supposed that the group
$\End_TM^{\times}$ posesses an $\O_T$-integral structure given by an appropriate
choice of an $\O_T$-lattice).}

Consider next the
${\bf C}_{p}$-analytic Lie group
$$
{\cal X}_{M,p} = {\rm Hom}_{contin}(GH_M,\     {\bf C}_{p}^{\times})
$$
 consisting of all continuous characters of the Hida group $GH_M$, 
which contains the ${\bf C}_{p}$-analytic Lie group
\par
$$
{\cal X}_{p} = {\rm Hom}_{contin}(\Gal_{p},\     {\bf C}_{p}^{\times})
$$
 consisting of all continuous characters of the Galois group $\Gal_{p}$
(via the projection of $GH_M$ onto $\Gal_p$.

The group ${\cal X}_{M,p}$ contains the discrete subgroup $\Ac$ of arithmetical 
characters of the type
$$
\chi \cdot \eta \cdot x_p^m = (\chi, \eta, m),
$$
where 
$$\chi \in {\cal X}_{M,p}^{tors}$$ is a character of finite order of $GH_M$,
 $\eta $ is a $T$-algebraic character of $(\End_TM)^{\times}(\O_{T, p})$, $m
\in \Z$, and $ x_p$ denotes the following natural homomorphism
$$
x_{p}: \Gal_{p}= \Gal(\Q_{p,\infty}^{ab}/{\bf Q}) \cong  {\bf
Z}^{\times }_{p} \rightarrow  {\bf C}^{\times }_{p},\     x_{p} \in  {\cal
X}_{p}. 
$$

\medskip
\noindent
{\bf Definition.} {\sl The {\it algebra of Iwasawa--Hida} $I_M = I_{M,p}$ of 
$M$ at $p$ is the completed group ring $\O_p[[GH_M]]$, where $\O_p$ denotes the
ring of integers of the Tate field $\C_p$.}

Note that this definition is completely analogous to the usual definition of 
the Iwasawa algebra $\La $ as the completed group ring $\Z_p[[\Z_p]]$ if we take
into account that $\Z_p$ coincides with the factor group of $\Z_p^{\times}$
modulo its torsion subgroup.

Now for each arithmetic point $P = (\chi, \eta, m) \in \Ac$ we have a 
homomorphism 
$$
\nu_P: I_{M,p} \rightarrow \O_p$$ which is defined by the
corresponding group homomorphism 
$$
P: GH_M \rightarrow \O_p^{\times} \subset
\C_p^{\times}.
$$

For a $I_M$-module $N$ and $P \in \Ac$ we put 
$$N_P = N \otimes _{I_M, \nu_P}\O_p$$ ("reduction of $N$ modulo $P$", or a
fiber of $N$ at $P$).  

Therefore, for a Galois representation
$$
r_N : G_\Q \rightarrow \GL(N)
$$
 of the above type
 its reduction $r_{N_P} = r \bmod P$ is defined as the natural composition:
$$
G_\Q \rightarrow \GL (N) \rightarrow \GL (N_P).
$$

\medskip
\noindent
{\bf Remark. } In his very recent work [HiGen] Hida gives another version of
the above definition, but he starts from a Galois representation 
$\ph :\Gal (\oe F/F) \longrightarrow GL_n({\bf I})$, where 
${\bf I}=\O_K[[T_n(\Z_p]]$ and $T_n$ the maximal split torus of $\Res
_{\O_F/\Z}GL(n)$ for the integer ring $\O_F$ of $F$, and for the integer 
ring $\O_K$ of a suffficiently large finite extension $K$ of $\Q_p$. 
He is interested in representations $\ph$ satisfying the following
condition: 

\medskip
{\it There are arithmetic points $P$ "densly populated" in $\Spec({\bf I}(K)$
such that the Galois 

representation $\ph_P=P\circ \ph$ is the $p$-adic \'etale
realization of a rank $n$ 
pure motive 

$M_P$ of weight $w$ defined over $F$ with
coefficients in a number field $E_P$ in $\oe Q$.}

\medskip
\noindent
We are trying to resolve an inverse problem and to include a given
motive $M$ in a maximal possible $p$-adic family $M_P$ parametrized by
arithmetic characters of a certain group which we suppose to consist of an
"algebraic part"  $(\End_TM)^{\times}(\O_{T, p})$ and 
of a "Galois part" $\Gal_p$.

\medskip
\noindent
{\bf  A conjecture on the existence of $p$-adic families of Galois
representations attached to motives. } 
Note first that the fixed embeddings $T
\hookrightarrow
\C$,                                                                                                                                                                                                    
$$  i_{\infty }: \oe \Q  \rightarrow  {\bf C},\ \    i_{p}: \oe \Q  \rightarrow 
{\bf C}_{p}
$$ define a place $\la (p)$ of $T$ attached to the corresponding
composition 
$$ T \hookrightarrow \oe \Q  \buildrel{i_p}\over\rightarrow  {\bf
C}_{p}. 
$$

\medskip
\noindent
{\bf Conjecture I. } {\sl For every $M$  of rang $d$ with coefficients 
in $T$ there esists a free $I_M$-module $M_I$ of the same rang $d$, a Galois
representation  $$r_I : G_F \rightarrow \GL(M_I),$$ a dense subset
$\Ac^\prime \subset \Ac$ of characters, and a distinguished point
$P_0 \in \Ac$ such that 

\par\noindent
(a) the reduced Galois representation 
$$r_{I,P_0} : G_F \rightarrow \GL(M_{I,P_0})$$
is equivalent over $\C_p$ to the $\la (p)$-adic representation $r_{M, \la (p)}$ of $M$ at the distinguished place $\la (p)$;

\par\noindent
(b) for every $P \in \Ac^\prime$ there exists a motive $M_P$ over $\Q$ of the
same rang $d$ such that its $\la (p)$-adic Galois representation is equivalent
over $\C_p$ to the reduction } 
$$ 
r_{I,P} : G_\Q\rightarrow \GL(M_{I,P}).
$$

We call the module $M_I$ the {\it realization of Iwasawa} of $M$.

\medskip
\noindent
{\bf  A generalization of the Hasse invariant for a motive. }
 We define the generalized {\it Hasse invariant} of a motive in terms  
of the {\it Newton polygons} and  the {\it
Hodge polygons} of a motive. Properties of these polygons are
closely related to the notions of a $p$-ordinary and a $p$-admissible
motive.

Now we are going to define the Newton polygon $P_{Newton}(u) =
P_{Newton}(u,  M)$ and the Hodge polygon 
$P_{Hodge }(u) = P_{Hodge }(u,  M)$ attached to $M$. 
First we consider (using $i_{\infty }$) the local ${
p}$-polynomial 
$$
L_{{ p}}(M ,  X)^{-1} = 1+A_{1}({ p})X + \cdots + A_{d}({ p})X^{d}
$$
$$
= (1-\alpha ^{(1)}({p})X)\cdot (1-\alpha ^{(2)}({ p})X)\cdot \dots \cdot
(1-\alpha ^{(d)}({ p})X),
$$
 and we assume that its inverse roots are indexed in such a way
that
$$
\hbox{ord}_{p}\alpha ^{(1)}({ p}) \le \hbox{ ord}_{p}\alpha ^{(2)}({p})\le
\cdots \le \hbox{ ord}_{p}\alpha ^{(d)}({ p}) 
$$

\medskip
\noindent
{\bf Definition. } {\sl The {\it Newton polygon} $P_{Newton}(u) (0\le  u \le d)$
of $M$ at ${p}$ is the convex hull of the points $(i,  
\ord_{p}A_{i}({p}))$ ($i = 0,  1,  \cdots ,   d$).}

The important property of the Newton polygon is that the length
 the horizontal segment of slope $i \in \Q$ is equal to the number of the
inverse roots $\alpha ^{(j)}({ p})$ such that ord$_{p}\alpha ^{(j)}({p}) = i
$ (note that $i$ may not necessarily be integer but this will be the case for
the $p$-ordinary motives below). \par
 The {\sl Hodge polygon} $P_{Hodge  }(u)$ $(0\le  u \le d)$ of $M$ 
is defined using the Hodge decomposition of the $d$-dimensional ${\bf
C}$-vector space $$
M_{B} = M_{B}\otimes _{T }{\bf C} = \oplus_{i, j}M^{i ,  j}
$$
 where  $M^{i ,  j}$ as a
${\bf C}$-subspace.  

\medskip
\noindent
 {\bf Definition. }{\sl  The {\it Hodge
polygon} $P_{Hodge}(u)$ is a function $[0, d] \rightarrow \R$ whose graph
consists of segments passing through the points  \par
$$
(0 ,  0),\     \dots ,\     (\sum_{i^\prime \le i}h(i^\prime  ,  j),\    
\sum_{i^\prime \le i} i^\prime h (i^\prime  ,  j)), 
$$
 so that the length of the horizontal segment of the slope $i \in \Z$ is
 equal to the dimension $h (i,  j)$.}
\par
Now we recall the definition of a $p$-ordinary motive (see [Co], [Co--PeRi]). 
We assume that $M$ is
pure of weight $w$ and of rank $d$. 
Let $G_p$ be the
decomposition group (of the place $\la (p)$ in $T$  over $p$) and
\par
$$
\psi _{p} :   G_p \rightarrow  {\bf Z}^{\times }_{p}
$$
 be the cyclotomic character of $G_p $. 
Then $M$ is called $p$ ordinary at
$p$ if the following conditions are satisfied:
\par\noindent
(i)  The inertia group $I_{p}\subset G_p $ acts trivially on each of the
$l$-adic realizations $M_{l}$ for $l\neq p$;
\par\noindent
(ii) There exists a decreasing filtration $F^{i}_{p} V$ on 
$V = M_{p} = M_{B}\otimes {\bf Q}_{p}$
of ${\bf Q}_{p}$-subspaces which are stable under the action of $G_p $ such that
for all $i\in {\bf Z}$ the group $G_p $ acts on $F^{i}_{p} V/F^{i+1}_{p}V$ 
via some power of the
cyclotomic character, say $\psi ^{-e_{i}}_{p}$. 
Then
\par
$$
e_{1}(M) \ge  \cdots \ge  e_{t}(M)
$$
 and the following properties take place:

\medskip\par\noindent
(a)
$$
 \dim _{{\bf Q}_{p}}F^{i}_{p} V/F^{i+1}_{p}V = h(e_{i} ,  w-e_{i});
$$

\medskip\noindent
(b) The Hodge polygon and the Newton polygon of $M$ coincide:
$$
P_{Newton}(u) = P_{Hodge}(u).
$$
If furthermore $M$ is critical at $s=0$ then it is easy to verify
that the number $d_{p}$ of the inverse roots $\alpha ^{(j)}(p)$ with
$$
\hbox{ord}_{p}\alpha ^{(j)}(p) < 0\hbox{ is equal to }d^{+}= 
d^{+}(M)\hbox{ of }M^{+}_B.
$$

   However, it turns out that the notion of a $p$-ordinary motive is 
too restrictive, and we have introduced the following weaker version of it.

\medskip
\noindent 
{\bf Definition. } {\sl The motive $M$ over $F$ with coefficients in $T$ 
is called
{\it admissible} at $p$ if
$$
P_{Newton }(d^{+}) = P_{Hodge}(d^{+})
$$
 here $d^{+}= d^{+}(M)$ is the dimension of the subspace $M^+\subset M_B$.
}
In the general case we use the following quantity (a "generalized slope")
$h = h_p$ which is defined as the difference
between the Newton polygon and the Hodge polygon of $M$:
$$
h_p = P_{Newton }(d^{+ }) - P_{Hodge}(d^{+ }). 
$$
 of $M$ at $p$.
Note the following important properties of $h$:
\par\noindent
(i) $h = h(M)$ does not change if we replace $M$ by its Tate twist.
\par\noindent
(ii) $h = h(M)$ does not change if we replace $M$ by its twist $M = M(\chi )$
with
a Dirichlet character $\chi $ of finite order whose conductor is prime to $p$.
\par\noindent
(iii) $h = h(M)$ does not change if we replace $M$ by its dual $M^{\vee}$
 \par
In the next section we state in terms of this quantity a
general conjecture on $p$-adic $L$-functions.

\medskip
\noindent
{\bf A conjecture on the existence of certain families of 
$p$-adic $L$-functions. }
We are going to describe families of $p$-adic $L$-functions as 
certain analytic functions on the total analytic space, the
${\bf C}_{p}$-analytic Lie group
$$
{\cal X}_{M,p} = {\rm Hom}_{contin}(GH_M,\     {\bf C}_{p}^{\times}),
$$
which contain the ${\bf C}_{p}$-analytic Lie subgroup (the cyclotomic line) 
${\cal X}_{p} \subset {\cal X}_{M,p}$:
$$
{\cal X}_{p} = {\rm Hom}_{contin}(\Gal_{p},\     {\bf C}_{p}^{\times}).
$$
In order to do this we need a modified $L$-function of a motive. 
Following J. Coates this  modified $L$-function has a form appropriate for
further use in the $p$-adic construction. First we multiply $L(M,  s)$ by an
appropriate factor at infinity and define 
$$ 
\Lambda _{(\infty )}(M,\     s) =
E_{\infty }(M,\     s)L(M ,  s) 
$$
 where $E_{\infty }(M,   s) = E_{\infty }(\tau,  R_{F/{\bf Q}}M,   \rho ,   s)$ 
is the modified $\Gamma $-factor at
infinity which actually does not depend on the fixed embedding $\tau$ of 
$T$ into $\C$. 
Also we put
$$
c ^{\nu }(M) = (c ^{\nu }(M)^{(\tau)}) = c^{\nu }(RM)(2\pi i)^{r(RM)}\in 
(T\otimes {\bf C})^{\times } 
$$
 where 
$$
\nu  = (-1)^{m}, r(M) = \sum_{j<0}jh(i,  j )
=\sum_{j<0}jh(i,  j), 
$$ 
$c^{\nu }(M)$ is the period of $M$. 
Note that the quantity $r(M)$ has a natural geometric interpretation as the 
minimum of the Hodge polygon $P_{Hodge}(M)$.

  We define 
$$
\eqalign{&
\Lambda _{(p, \infty )}(M(m)(\chi),     s) = \cr &
 G(\chi )^{-d^{\ep _{0}}(M(m)(\chi))}\prod_{{ p} \mid p} A_{{
p}}(M(m)(\chi), s)\cdot \Lambda _{(\infty )}(M(m)(\chi),    s)}, 
$$
where
$$
A_{p}(M(\chi ),s) = \cases{\prod_{i = d^+ +1}^d(1 - \chi (p )\al
^{(i)}(p) p^{-s})\prod_{i = 1}^{d^+}(1 - \chi ^{-1}(p )\al
^{(i)}(p)^{-1} p^{s-1})&\cr &\cr
 \hskip3.5cm{\rm \ \ for}\  p \nd \c(\chi)&\cr &\cr \prod_{i = 1}^{d^+}
{\displaystyle\left(\frac p^s,\al ^{(i)}(p),\right)^{\ord_{p}\c(\chi )}},
{\rm otherwise.}&}
$$

\medskip
Let $\Ac$ be the discrete subgroup $\Ac$ of arithmetical characters,
$$
\chi \cdot \eta \cdot  x_p^m = (\chi, \eta, m) \in \Ac,
$$
$\Ac^\prime \subset \Ac$ a certain "dense" subset of characters,   $P_0
\in \Ac$ a distinguished point of conjecture I. 
Let $\Ac^{\prime
\prime} \subset  \Ac^\prime$ be the subset of critical elements, which consists
of those $P$, for which the corresponding motives $M_P$ are critical (at $s =
0$). 
Now we are ready to formulate the following

\medskip
\noindent
{\bf Conjecture II. } {\sl 
There exists a certain choice of complex periods  $\W _{\infty}(P) \in
\C^{\times}$ and $p$-adic periods $\W_p(P) \in \C_p^{\times}$ for all $P
\in \Ac^{\prime \prime}$ such that ``the ratio''  $\W_{p}(P)/\W_\infty (P)$ is
canonicaly defined, and there exists a $\C_p$-meromorphic  function 
$$ 
\L_M : {\cal X}_{M,p} \rightarrow \C_p
$$ with
the properties:
\medskip\par\noindent
(i) 
$$
 \L_M (P) = \W_p(M) {\Lambda _{p, (\infty )}(M(m)(\chi),     0)\over
\W_\infty (P)} 
$$ for almost all $P \in \Ac^{\prime \prime}$;
\medskip\par\noindent
(ii) For  arithmetic points of type 
$$
P = (\chi , \eta, m)\in \Ac^{\prime \prime}
$$ 
with $\eta$ fixed there exists a finite set $\Xi  \subset  {\cal X}_{M,p}$ of 
$p$-adic characters and
positive integers $n(\xi ) $ (for $\xi \in \Xi )$ such that for any 
$g_{0} \in \Gal _{p}$ we
have that the function 
$$
\prod_{\xi \in \Xi } (x(g_{0})- \xi (g_{0}))^{n(\xi )} \L_{M}(x\cdot P)
$$
 is holomorphic on ${\cal X}_{p}$;
\medskip\par\noindent
(iii) For  arithmetic points of type 
$$
P = (\chi , \eta, m)\in 
\Ac^{\prime \prime}
$$ with $\eta$ fixed the function  in (ii) is bounded if and
only if the invariant $h (P) = h(M_P)$ vanishes;
\medskip\par\noindent
(iv) In the general case the function $\L_M(P\cdot x)$ of $x \in \X_p$
is of logarithmic growth type $o(\log  (\cdot )^h{_0})$ with 
$$
h_0 = [ h] + 1.
$$ 
\medskip\par\noindent
(v) For  arithmetic points of type 
$$
P = (\chi , \eta, m)\in 
\Ac^{\prime \prime}
$$ with $\chi$ and $m$ fixed the function  in (ii) is always bounded if the
Hasse invariant $h (P) = h(M_P)$ does not depend on $\eta$. 
} 
 
\medskip
\noindent
Note that the assertion (v) means in particulary that
 the values of the function 
$$
P\mapsto \W_p(M) {\Lambda _{(p, \infty )}(M_P,     0)\over
\W_\infty (P)} 
$$ 
 satisfy 
 generalised Kummer congruences in the followinfg sense:
 for any finite linear combination $\sum_Pb_P\cdot P$ with $b_P\in \C_p$ which
 has the property 
 $\sum_Pb_P\cdot P\equiv 0 (\bmod p^N)$ we have that for some 
 constant $C\not = 0$ the
 corresponding
 linear combination of the normalized $L$-values
$$
C\sum _P b_P \W_p(M_P) \cdot  {\Lambda _{(p, \infty )}(M_P,    
0)\over \W_\infty (P)}  
\equiv 0 (\bmod p^N). 
$$

\medskip\noindent In  the case of families of supersingular modular forms
studied by R.Coleman [CoPBa] the invariant $h(P)$ reduces to the slope of a
modular form in such a family.

\bigskip
\bigskip
\noindent
{\bf References}
\medskip

\par
\medskip
\newdimen\brawidth\setbox0=\hbox{[PaAdm1]\ }\brawidth=\wd0

\def\ref#1,#2!{\smallskip\noindent
\hangindent=\brawidth
\hangafter=1
\hbox to\brawidth{[#1]\hfill}{#2}:\ }

\ref Am--V,Amice Y.,Velu J.!Distributions $p$-adiques associ\'ees aux
s\'eries de Hecke, Jour\-n\'ees Arithm\'etiques de Bordeaux (Conf. Univ.
Bordeaux, 1974), Ast\'e\-risque no. 24/25, Soc. Math. France, Paris
1975, 119--131.

\ref Co,Coates J.!On $p$-adic $L$-functions. S\'em. Bourbaki, 40\`eme ann\'ee,
1987-88, no 701, Asterisque (1989) 177--178.

\ref Co--PeRi,Coates J., Perrin-Riou B.!On $p$-adic $L$-functions
attached to motives over ${\bf Q}$. Advanced Studies in Pure Math. 17
(1989), 23--54.

\ref CoPBa,Coleman R.!$P$-adic Banach Spaces and Families of Modular forms.
Manu\-script of January 7, 1995.

\ref GouMa,Gouv\^ea F.Q. and Mazur B.!On the characteristic power series of the
$U$-operator, Ann. Inst. Fourier, Grenoble 43, 2 (1993), 301--312.

\ref  HiGal,Hida H.!Galois representations into $GL_{2}({\bf Z}_{p}[[X]])$ 
 attached to ordinary cusp forms,
 Invent. Math. 85 (1986) 545--613.

\ref HiGen,Hida H.!On the search of genuine $p$-adic modular $L$-functions 
 for $\GL (2)$. Manuscript of December 12, 1994.

\ref KapCM,Katz N.M.!$p$-adic $L$-functions for $CM $-fields, Invent.
 Math. 48 (1978) 199--297.
 
\ref Maz,Mazur B.!Deforming Galois representations // Galois Groups over $\Q$.
Ed. Y.Ihara, K.Ribet, J.-P. Serre, 1989, Springer.

\ref Maz--W,Mazur B., Wiles A.!On $p$-adic analytic families of Galois
 representations // Compos. Math. 59 (1986) 231--264.

\ref PaLNM,Panchishkin A.A.!Non-Archimedean $L$-functions of Siegel
 and Hilbert modular forms, Lecture Notes in Math., 1471, 
 Springer, 1991.

\ref PaAdm,Panchishkin A.A.!Admissible Non-Archimedean standard zeta functions
 of Siegel modular forms, 
 Proceedings of the Joint AMS Summer Conference on Motives,
 Seattle 1991, AMS, Providence, R.I., 1993, vol.2, 
 251--292.

\ref PaIF,Panchishkin A.A.!Motives over totally real fields and $p$-adic
 $L$-functions,  Annales de
 l'Institut Fourier, Grenoble, 44, 4 (1994).

\vskip2cm

{\font\caps=cmcsc9\parindent=0pt

\caps
Mathematical Science Research Institute, Berkeley, CA 94720, USA

\medskip
{\rm and}
\medskip

Institut Fourier, B.P.74, 38402 St.-Martin d'H\`eres, France}

\end